\documentclass[10pt, a4paper, oneside, reqno]{amsart}

\makeatletter
\def\@settitle{\begin{center}%
		\baselineskip14\p@\relax
		\normalfont\LARGE\bfseries
		\@title
	\end{center}%
}

\def\section{\@startsection{section}{1}%
	\z@{.7\linespacing\@plus\linespacing}{.5\linespacing}%
	{\normalfont\large\bfseries}}

\def\subsection{\@startsection{subsection}{2}%
	\z@{.5\linespacing\@plus.7\linespacing}{.5\linespacing}%
	{\normalfont\bfseries}}

\def\@setauthors{%
  \begingroup
  \def\thanks{\protect\thanks@warning}%
  \trivlist
  \centering\footnotesize \@topsep30\p@\relax
  \advance\@topsep by -\baselineskip
  \item\relax
  \author@andify\authors
  \def\\{\protect\linebreak}%
  \authors%
  \ifx\@empty\contribs
  \else
    ,\penalty-3 \space \@setcontribs
    \@closetoccontribs
  \fi
  \endtrivlist
  \endgroup
}

\makeatother

\usepackage[table, xcdraw, usenames, dvipsnames]{xcolor}
\definecolor{darkblue}{rgb}{0.0, 0.0, 0.45}
\definecolor{darkgreen}{rgb}{0.0, 0.45, 0}
\usepackage[colorlinks	= true,
raiselinks	= true,
linkcolor	= darkblue, 
citecolor	= Mahogany,
urlcolor	= darkgreen,
pdfauthor	= {},
pdftitle	= {},
pdfkeywords	= {},
pdfsubject	= {},
plainpages	= false]{hyperref}

\pdfoutput=1
\date{\today}

\usepackage{dsfont,amsfonts,amssymb,amsmath,amsthm}
\usepackage{mathtools, mathrsfs}

\usepackage{graphicx}
\usepackage[font=small,labelfont=rm]{subcaption}
\usepackage[font=small,margin=10pt]{caption}
\usepackage{wrapfig}

\usepackage{multirow}
\usepackage{booktabs}

\usepackage{algorithm,algorithmic}

\usepackage{enumitem}
\usepackage{tikz}
\usepackage{courier} 

\usepackage{nicefrac}

\theoremstyle{plain}
\newtheorem{Thm}{Theorem}[section]

\newtheorem{Prop}[Thm]{Proposition}

\newtheorem{Lem}[Thm]{Lemma}

\newtheorem{As}[Thm]{Assumption}

\newtheorem{Rem}[Thm]{Remark}

\DeclareMathOperator*{\argmin}{\arg\!\min}
\DeclareMathOperator*{\argmax}{\arg\!\max}

\newcommand{\R}{\mathbb{R}}
\newcommand{\N}{\mathbb{N}}

\newcommand{\Rl}{\underline {\R}}

\newcommand{\ra}{\rightarrow}

\newcommand{\wh}{\widehat}
\newcommand{\Let}{\coloneqq}

\newcommand{\tr}{^{\top}}

\newcommand{\norm}[1]{\left\Vert #1 \right\Vert}

\newcommand{\EE}{\mathds{E}}

\DeclareMathOperator{\ord}{\mathcal{O}}
\newcommand{\ind}[1]{\mathds{1}_{#1}}

\newcommand{\opt}{_\star}
\newcommand{\X}{\mathsf{X}}
\newcommand{\Y}{\mathsf{Y}}
\newcommand{\U}{\mathsf{U}}
\newcommand{\Z}{\mathsf{Z}}
\newcommand{\W}{\mathsf{W}}

\newcommand{\sam}{_{\mathrm{s}}}
\newcommand{\dpo}{\mathbf{B}}
\newcommand{\sdp}{\mathbf{B}\sam}

\newcommand{\lf}{L}
\newcommand{\rew}{r}

\usepackage[square,numbers]{natbib}

\graphicspath{{Fig/}}

\title[]{Fitted Q-Iteration via Max-Plus-Linear Approximation}

\author[]{Y. Liu and M.~A.~S.~Kolarijani}
\thanks{The authors are with Delft Center for Systems and Control, Delft University of Technology, 
Delft, The Netherlands. Email: \texttt{\{Y.Liu-12, M.A.SharifiKolarijani\}@tudelft.nl}.}
\thanks{This research was partially supported by the funding received from the European Research Council (ERC) under the grant TRUST-949796.}

\begin{document}

\begin{abstract}
In this study, we consider the application of max-plus-linear approximators for Q-function in offline reinforcement learning. 
In particular, we incorporate these approximators to propose novel fitted Q-iteration (FQI) algorithms with provable convergence. 
Exploiting the compatibility of the Bellman operator with max-plus operations, we show that the max-plus-linear regression within each iteration of the proposed FQI algorithm reduces to simple max-plus matrix-vector multiplications. 
We also consider the variational implementation of the proposed algorithm which leads to a per-iteration complexity that is independent of the number of samples. 

\smallskip

\end{abstract}

\maketitle

\section{Introduction}\label{sec:intor}

Reinforcement learning (RL) is a subfield of machine learning concerned with sequential decision-making in uncertain environments~\cite{Sutton18, Bertsekas19}. 
Over the past three decades, RL has gained popularity in various applications such as robotics~\cite{Levine16}, game playing~\cite{silver2017mastering}, and finance~\cite{Spooner18market}. These and other successful applications of RL showcase its ability to adapt and optimize behaviors in uncertain environments requiring a high degree of strategic depth. 

The ultimate goal of RL is for an agent (a.k.a. decision-maker, controller) to learn to perform sequences of actions that maximize the cumulative reward. 
In \emph{offline} RL (a.k.a. batch RL), the agent has to learn a (control) policy from a fixed dataset of interactions with the environment, as opposed to being able to interact with the environment in an online fashion. 
This approach is particularly useful when online interaction is expensive, risky, or infeasible, e.g., in healthcare applications~\cite{nie2021lea}. 

A standard approach for solving the offline RL problem is the fitted Q-iteration (FQI) algorithm~\cite{ernst05a}. 
FQI forms an approximation of the Q-function (i.e., state-action value function) and iteratively updates it by minimizing the empirical Bellman error using the available batch of data. 
The specific approximator used in FQI can be parametric (e.g., neural networks~\cite{rie2005neu}) or non-parametric (e.g., kernel-based~\cite{far2009reg}). 
While the simplicity of FQI makes it a popular choice in practice, this algorithm lacks guarantees for convergence. 
Indeed, FQI can diverge even for problems with a finite number of states and actions and with a linear parametric approximator that is rich enough to represent the optimal value function~\cite{BAIRD199530, tsitsiklis1996feature}. 
A sufficient but restrictive condition for the convergence of FQI is for the approximator to be an ``averager''~\cite{GORDON1995261}, which is for instance the case for piece-wise constant approximations.

A particular class of approximators with a long history in optimal control and dynamic programming (DP) problems is \emph{max-plus (MP)}-linear approximators. 
These approximators, as the name suggests, are linear in the MP algebra in which the conventional addition and multiplication operations are replaced by maximum and addition operations, respectively. 
This widespread application stems from the compatibility of the Bellman operator with MP operations. 
Indeed, the backward iteration in the DP algorithm for solving the optimal control problem of finite deterministic Markov decision processes (MDPs) is a linear operator in MP algebra \cite{gaubert1997methods}. 
Another example is the MP-linearity of the evolution semigroup of the Hamilton-Jacobi-Bellman equation, arising in optimal control of deterministic continuous-time dynamics, that has been exploited to derive efficient numerical methods for finding an MP-linear approximation of the corresponding value function~\cite{fle2000max, McEneaney03, Akian08}. 
More recently, a similar idea has been used for approximating the value function of a deterministic MDP with a continuous state space~\cite{Bach20}. 
Also of interest is the recent development of an \emph{online} RL algorithm that, unlike the previous works mentioned above, updates the parameters of the MP-linear approximation of the value function in real-time~\cite{goncalves2021maxplus}. 

In this paper, we propose a novel class of FQI algorithms that employ parametric MP-linear approximation of the Q-function. 
This, in turn, allows us to exploit the MP-linearity of the empirical Bellman operator to reduce the MP-linear regression problem in each iteration to a simple max-plus matrix-vector multiplications. The main \textbf{contributions} of this work are as follows: 
\textbf{(i)}~We propose the MP-FQI Algorithm~\ref{alg:mp-fqi} and show its convergence with a linear rate (Theorem~\ref{thm:mp-fqi conv}) and a per-iteration complexity of $\ord(np)$, where $n$ is the number of samples and $p$ is the number of parameters/basis functions in the MP-linear approximator (Theorem~\ref{thm:mp-fqi comp}). 
\textbf{(ii)}~Inspired by~\cite{Akian08}, we also consider the variational formulation of the problem and propose the v-MP-FQI Algorithm~\ref{alg:mp-fqi-var} and show its convergence with a linear rate (Theorem~\ref{thm:mp-vqi conv}) and a per-iteration complexity of $\ord(pq)$, independent of the number of samples, where $p$ is again the number of parameters and $q$ is the number of test functions in the variational formulation (Theorem~\ref{thm:mp-vqi comp}). 
To the best of our knowledge, this is the first work that provides a provably convergent FQI algorithm using a class of approximators beyond the existing restrictive conditions that require the mapping corresponding to the approximator to be non-expansive in the $\infty$-norm such as~\cite{GORDON1995261}.

The rest of the paper is organized as follows. In Section~\ref{sec:RL and FQI}, we provide the statement of the offline RL problem and its standard solution via FQI. 
We also provided some preliminaries on MP-linear approximation and regression. 
The MP-FQI algorithm and its variational version along with the convergence and complexity analysis are provided in Sections~\ref{sec:MP-FQI} and \ref{sec:variational}, respectively. 
In Section~\ref{sec:num ex}, we examine the performance of the proposed algorithms in comparison with the ``standard'' FQI algorithm through numerical examples. 
Section~\ref{sec:conclusion} concludes the paper with some final remarks. 
All the technical proofs are provided in Appendix~\ref{app:proof}. 

\textbf{Notations.} 
$\R$ denotes the set of real numbers and $\Rl = \R \cup \{ -\infty \}$ is the extended real line. 
We denote $[1:n] = \{ 1,\ldots, n \}$ for $n \in \N$, where $\N = \{1,2,3,\ldots\}$ is the set of natural numbers. We also denote  $\N_0 = \{0,1,2,3,\ldots\}$. 
Given a set $\X$, we define the indicator function~$\ind{\X}(x) = 1$ if $x\in\X$ and $=0$ otherwise, and the \emph{max-plus} indicator function~$\delta_{\X}(x) = 0$ if $x\in\X$ and $=-\infty$ otherwise. 
We use $[v]_i$ to denote the $i$-th element of the vector $v \in \Rl^n$, and $[M]_{i,j}$ to denote the entry on row~$i$ and column~$j$ of the matrix $M \in \Rl^{n\times m}$. 
$\norm{v}_2 = \sqrt{v\tr v}$ and $\norm{v}_{\infty} = \max_{i\in[1:n]} |[v]_i|$ denote the 2-norm and $\infty$-norm of the vector $v \in \R^n$, respectively. 
The all-one vector is denoted by $e = (1,1,\ldots, 1)\tr$ and the identity matrix is denoted by $I$. 
A matrix~$M \in \Rl^{n\times m}$ is called column (resp.~row) $\R$-astic  if it has at least one finite entry in each column (resp.~row). 
A matrix that is both column and row $\R$-astic is called doubly $\R$-astic. 
We use bold capital letters as in $\mathbf{A}, \mathbf{B}, \mathbf{D}$ to denote operators. 
To simplify the presentation, we particularly use $\boxtimes$ to denote the max-plus matrix multiplication, i.e., for two matrices $A \in \Rl^{n\times m}$ and $B \in \Rl^{m \times p}$, we define $A\boxtimes B \in \Rl^{n\times p}$ with entries $[A\boxtimes B]_{i,j} = \max_{k \in [1:m]} [A]_{i,k}+[B]_{k,j}$.

\section{Problem statement and preliminaries}
\label{sec:RL and FQI}

In this section, we provide the problem statement and its standard solution. 
We also review the background on MP-linear approximation and regression along with some preliminary lemmas to be used in the subsequent sections.

\subsection{Offline RL} 
Consider an MDP~$\big(\Z = (\X\times\U), P, r, \gamma\big)$ where $\X\subset \R^{d_{\mathrm{x}}}$ and $\U\subset \R^{d_{\mathrm{u}}}$ are the state and action spaces, respectively; $P:\Z\times\X\ra[0,\infty)$ is the transition probability kernel such that $P\big(z = (x,u), \cdot\big) $ is a probability measure on $\X$ describing the distribution of the next state~$x^+$ given that the action~$u$ is taken in state~$x$; $ r:\Z\ra\R$ is the reward function and uniformly bounded from above such that $r\big(z = (x,u)\big)$ is the reward that the agent receives for taking the action~$u$ in state~$x$; $\gamma \in (0,1)$ is the discount factor.

The problem of interest is to find a control policy that maximizes the expected, infinite-horizon, discounted reward starting from each stat-action pair, that is, to find the optimal Q-function $Q\opt : \Z \ra \R$ (a.k.a.~optimal state-action value function) given by
\begin{align*}
   Q\opt(z) \Let \max_{\pi:\X\ra\U}  \EE \left( \sum_{t=0}^{\infty} \gamma^t r(z_t)\ \big|\ z_0 = z,\  u_t = \pi(x_t),\ \forall t\geq 1 \right). 
\end{align*}
The optimal policy is then $\pi\opt(x) \in \argmax_{u\in\U} Q\opt(x,u)$, i.e., \emph{the greedy policy with respect to $Q\opt$}. 
The offline RL problem involves finding the Q-function using a set of $n\in\N$ samples 
$ \{ z_i=(x_i, u_i), x^+_i, r_i \}_{i\in [1:n]}$, 
where each sample corresponds to the agent taking the action~$u_i$ in the state~$x_i$, which pushes the system to the state~$x^+_i$ while receiving the reward $r_i$. 

\subsection{Fitted Q-iteration (FQI)}
The optimal Q-function can be characterized as the fixed-point of the Bellman operator~$\dpo:\R^{\Z}\ra \R^{\Z}$, 
that is, $Q\opt$ solves the Bellman equation (BE) \cite[Fact~3]{szepesvari10}:  
\begin{equation} \label{eq:DP}
Q(z) = \dpo Q(z) \Let \EE_{x^+} \big(  r(z) + \gamma \cdot \max_{u^+} Q(x^+,u^+)\big),\quad \forall z \in \Z.
\end{equation}
Based on this result, a standard solution for the offline RL problem is to use a parametric approximation $Q_{\theta}\in \R^{\Z}$ and then find the parameter $\bar{\theta} \in \R^p$ such that $Q_{\bar{\theta}}$ solves the \emph{empirical} BE:
\begin{equation} \label{eq:DP_empirical}
Q_{\theta}(z_i) = \sdp Q_{\theta}  (z_i) \Let  r_i + \gamma \cdot \max_{u^+} Q_{\theta} (x_i^+,u^+),\quad \forall i \in [1:n].
\end{equation}
(We use the subscript $\mathrm{s}$, for the sampled versions of the same objects). 
The so-called FQI algorithm solves the latter problem by minimizing the empirical Bellman error in a recursive fashion~\cite{ernst05a}. 
To be precise, with some initialization $\theta^{(0)}$, FQI solves the regressions
\begin{equation}\label{eq:fqi_stand}
\theta^{(\ell+1)} \in    \argmin_{\theta} \sum_{i\in [1:n]} \lf\big(Q_{\theta}(z_i), \sdp Q_{\theta^{(\ell)}}(z_i)\big),  \quad \ell\in\N_0, 
\end{equation}
where $L:\R\times\R\ra\R$ is a loss function.
That is, FQI finds $\theta^{(\ell+1)}$ by fitting $Q_{\theta^{(\ell+1)}}$ to $\sdp Q_{\theta^{(\ell)}}$ over the sampled state-action space $\Z\sam \Let \{z_i\}_{i\in[1:n]}$. 

Alternatively, we can think of the procedure explained above as finding the fixed-point of the composite operation $\mathbf{R}\circ \dpo \circ \mathbf{A}: \R^p \ra \R^p$ on the space of parameters, 
with $\mathbf{A}: \R^p \ra \R^{\Z}$ corresponding to the function approximation ($\theta^{(\ell)} \mapsto Q_{\theta^{(\ell)}}$) and 
$\mathbf{R}: \R^{\Z} \ra \R^p$ corresponding to the regression~\eqref{eq:fqi_stand} \cite[Sec.~3.4.4]{Busoniu10}. 
Consequently, the convergence of the algorithm depends on the properties of the operators $\mathbf{A}$ and $\mathbf{R}$. 
In particular, if these two operations are non-expansive (in ${\infty}$-norm), then the whole iterative process becomes convergent since $\dpo$ (and hence $\sdp$) is a $\gamma$-contraction (in ${\infty}$-norm). 
In the subsequent sections, we will use this basic procedure in order to develop similar algorithms using MP-linear approximation and regression.

\subsection{MP-linear approximation}  
For a function $Q \in \R^{\Z}$, the MP-linear approximator is of the form   
\begin{equation}\label{eq:mp-approx}
    Q_{\theta}(z) = \max_{j\in [1:p]} \left\{ [f]_j(z) + [\theta]_j \right\} = f(z)\tr \boxtimes \theta,
\end{equation}
where $f:\Z \ra \Rl^p$ is the vector of MP features (basis functions), 
and $\theta \in \R^p$ is the vector of MP parameters (coefficients). 
For the features, two of the classical choices are (i)~the \emph{MP indicator function}~$ [f]_j(z) = \delta_{\W_j} (z)$ with $\{\W_j\}_{j \in [1:p]}$ being a partitioning of the set $\Z$ (leading to a \emph{piece-wise constant} approximation) and (ii)~\emph{quadratic} functions~$[f]_j(z) = -c \norm{z-w_j}_2^2$ with $c>0$ and $\{w_j\}_{j \in [1:p]} \subset \Z$ (leading to a \emph{semi-convex, piece-wise quadratic} approximation~\cite[Sec.~2]{McEneaney06}). 
In~\cite{Bach20}, the authors combine these two functions and use the \emph{distance} features~$[f]_j(z) = -c \min_{w\in \W_j} \norm{z - w}_2^2$ leading to an \emph{almost} piece-wise constant approximation.

\subsection{MP-linear regression}  
Consider the system of MP-linear equations $A\boxtimes \theta = b$ in the unknown $\theta \in \Rl^p$, with data $A \in \Rl^{n\times p}$ and $b \in \R^n$ such that the matrix $A$ is doubly $\R$-astic, i.e., it has at least one finite entry in each column and one finite entry in each row~\cite{cuninghame79}.  
This equation has a solution if and only if the \emph{principal solution} 
\begin{equation}
    \tilde{\theta} \Let - \left( A\tr \boxtimes (-b) \right) \in \R^{p},
\end{equation}
satisfies $A\boxtimes \tilde{\theta} = b$~\cite[Cor.~3.1.2]{butkovivc10}. 
Let $\varTheta = \{\theta \in \Rl^p: A\boxtimes \theta \leq b \}$ be the set of all the points that satisfy the so-called ``lateness'' constraint~\cite{cuninghame79}. 
The principal solution is also the \emph{greatest subsolution} in the sense that $\tilde{\theta} \in \varTheta$ and $\tilde{\theta} \geq \theta$ for all $\theta \in \varTheta$~\cite[Thm.~3.1.1]{butkovivc10}. 
Observe that the operation corresponding to the principal solution, i.e., the mapping $b\mapsto \tilde{\theta} = - \left( A\tr \boxtimes (-b) \right)$, can be equivalently seen as a projection.
Moreover, the principal solution is the optimal solution of the constrained $\infty$-norm regression~\cite[Thm.~3.5.1]{butkovivc10}, that is, 
$$
\tilde{\theta} = \argmin_{\theta\in \Rl^p} \left\{ \norm{A\boxtimes \theta - b}_{\infty}: A\boxtimes \theta \leq b \right\}.
$$
Another strong result concerns the optimal solution of the unconstrained $\infty$-norm regression as follows
$$\wh{\theta} = \argmin_{\theta\in \Rl^p} \norm{A\boxtimes \theta - b}_{\infty} = \tilde{\theta} + \frac{\Vert A\boxtimes \tilde{\theta} - b\Vert_{\infty}}{2}  e,$$ 
where $e$ is the all-one vector~\cite[Thm.~3.5.2]{butkovivc10}. 

\subsection{Preliminary lemmas}
We finish this section with two preliminary results. 
In the algorithms to be developed in the subsequent sections, we will exploit the following two important properties of the empirical Bellman operator. 
We note that the following result has been extensively studied and used in previous works such as \cite{gaubert1997methods, fle2000max, McEneaney03, Akian08, Bach20, goncalves2021maxplus}.  

\begin{Lem}[MP additivity and homogeneity of $\sdp$] \label{lem:mp add and hom} 
Consider the two functions $f,\tilde{f} \in \Rl^{\Z}$ and the scalar $\alpha \in \Rl$. 
Define $[\max \{ f, \tilde{f} \}](z) = \max\{f(z),\tilde{f}(z)\}$ and $[\alpha + f](z) = \alpha + f(z)$ for all $z\in\Z$. 
We have 
\begin{align}
&\sdp [\max \{ f, \tilde{f} \}]  = \max \{ \sdp f , \sdp \tilde{f} \}, \tag{MP additivity} \\
&\sdp [\alpha + f] = \gamma \alpha + \sdp f \tag{MP homogeneity}.
\end{align}
\end{Lem}

The next result concerns a generalization of the non-expansiveness of MP-linear operators. 
We will use this result for the convergence analysis of the algorithms proposed in this study. 
 
\begin{Lem}[Non-expansiveness of MP-linear operators]\label{lem:contraction general}
Consider the operator $\mathbf{A}: \R^n \ra \R^m: y \mapsto A \boxtimes (\eta y)$, where $A \in \Rl^{m \times n}$ is row $\R$-astic and $\eta \in (0,1]$ is a constant.  Then, 
\begin{equation}
\norm{\mathbf{A}y_1 - \mathbf{A}y_2}_{\infty} \leq \eta  \norm{y_1 - y_2}_{\infty}, \quad \forall y_1,y_2 \in \R^n.
\end{equation} 
\end{Lem}

\section{Max-plus FQI}
\label{sec:MP-FQI}

In this section, we propose and analyze the max-plus fitted Q-iteration (MP-FQI) algorithm. 
The basic idea here is to approximate the Q-function as an MP-linear function and use  MP-linear regression to find the parameters. 

\subsection{Algorithm} 
Let us consider an MP-linear approximation of the Q-function given by  
\begin{equation}\label{eq:mp-approx Q}
    Q_{\theta}(z) =  \max_{j\in [1:p]} \left\{ [f]_j(z) + [\theta]_j \right\} = f(z)\tr \boxtimes \theta,
\end{equation}
with features $f:\Z \ra \Rl^p$ and parameters $\theta \in \R^p$. 
Recall that the FQI algorithm solves the empirical BE
\begin{equation}\label{eq:DP_empirical_2}
Q_{\theta}(z_i) = \sdp Q_{\theta}  (z_i),\quad \forall i \in [1:n],
\end{equation}
via recursive regressions. 
The following lemma provides an alternative characterization of the preceding equation.

\begin{Prop}[MP empirical BE]\label{prop:MP_Bellman}
Define the matrices $F\sam , G\sam \in  \Rl^{n\times p}$ with entries
\begin{equation*}
    [F\sam]_{i,j} = [f]_j(z_i), \  [G\sam]_{i,j} = \sdp [f]_j(z_i), \qquad i\in[1:n],\ j\in[1:p], 
\end{equation*}
and let $\bar{\theta} \in \R^p$ be a solution to the MP equation
\begin{equation}\label{eq:DP_empirical_mp}
    {F\sam} \boxtimes \theta = G\sam \boxtimes (\gamma\theta).
\end{equation}
Then, the MP-linear function $Q_{\bar{\theta}}(z) = f(z)\tr \boxtimes \bar{\theta}$ solves the empirical BE~\eqref{eq:DP_empirical_2}. 
\end{Prop}

The preceding result implies that finding the parameters of the MP-linear approximation~\eqref{eq:mp-approx Q} of the Q-function reduces to solving the MP equation~\eqref{eq:DP_empirical_mp}. 
Similar to the standard FQI, we use recursive MP-linear regressions of the form 
\begin{equation}\label{eq:MP-FQI_iter0}
    \theta^{(\ell+1)} = \argmin_{\theta}\Vert{F\sam} \boxtimes \theta -  G\sam \boxtimes (\gamma\theta^{(\ell)})\Vert_{\infty}  , \quad \ell \in \N_0.
\end{equation}
To be able to solve these regressions efficiently, we need the following assumption:

\begin{As}
\label{As:features} 
The features $f:\Z\ra\Rl^p$ and the samples $ \{ z_i=(x_i, u_i), x^+_i, r_i \}_{i \in [1:n]}$ are such that 
\begin{align}
    &\forall z \in \Z,\ \exists j \in [1:p],\ [f]_j(z) \neq -\infty; \label{eq:F proper} \\
    &\forall j \in [1:p],\ \exists i \in [1:n],\ [f]_j(z_i) \neq -\infty. \label{eq:F col realistic}
\end{align}
That is, $F\sam$ (resp., $G\sam$) is doubly (resp., row) $\R$-astic.
\end{As}

Note that the preceding assumption is only relevant when the features are \emph{extended} real-valued functions that assign  $-\infty$ to a subset of state-action space~$\Z = \X\times\U$. 
In particular, if condition~\eqref{eq:F proper} does not hold, then any MP-linear combination of features assigns $-\infty$ to a subset of $\Z$. 
On the other hand, condition~\eqref{eq:F col realistic} implies that each feature $[f]_j$ is ``activated'' by at least one state-action sample $z_i$ in the sense that $[f]_j(z_i)$ is finite. 
If this condition does not hold for some $j \in [1:p]$, we can safely remove the corresponding feature~$[f]_{j}$ since the sample set is not ``rich'' enough to activate this feature (the algorithm below sets $[\theta]_j = -\infty$ in all iterations). 
\begin{Prop}[MP-FQI regression]\label{prop:MP_FQI_reg}
Let Assumption~\ref{As:features} hold. Then, the solution to the MP-linear regression~\eqref{eq:MP-FQI_iter0} is given by
\begin{equation*}
    \theta^{(\ell+1)} = \tilde{\theta}^{(\ell+1)} + \beta^{(\ell+1)} e  , \quad \ell \in \N_0,
\end{equation*}
where
\begin{align}
    &\tilde{\theta}^{(\ell+1)} = \mathbf{D}\tilde{\theta}^{(\ell)} \Let  - \left[ F\sam\tr \boxtimes \big[-(G\sam \boxtimes (\gamma\tilde{\theta}^{(\ell)}))\big] \right], \label{eq:MP-FQI_iteration} \\
    &\beta^{(\ell+1)} = \gamma  \beta^{(\ell)} + \frac{1}{2} \Vert{F\sam} \boxtimes \tilde{\theta}^{(\ell+1)} -  G\sam \boxtimes (\gamma\tilde{\theta}^{(\ell)})\Vert_{\infty}, \label{eq:MP-FQI_const}
\end{align}
with $\tilde{\theta}^{(0)} = \theta^{(0)}$ (i.e., $\beta^{(0)} = 0$).  
\end{Prop}
Now, observe that because of the MP-linear structure of the approximation $Q_{\theta^{(\ell+1)}}$, any constant shift such as $\beta^{(\ell+1)}$ in the entries of $\theta^{(\ell+1)}$ leads to the same greedy policy; 
that is, the greedy policies with respect to $Q_{\theta^{(\ell+1)}}$ and $Q_{\tilde{\theta}^{(\ell+1)}}$ are the same. 
Therefore, it suffices to follow the iterations of principal solution~$\tilde{\theta}$ in~\eqref{eq:MP-FQI_iteration}. 
Algorithm~\ref{alg:mp-fqi} summarizes the procedure described above.

\begin{algorithm}[t]
\begin{small}
   \caption{Max-plus FQI (MP-FQI)} 
   \label{alg:mp-fqi}
\begin{algorithmic}[1]
	\REQUIRE samples~$ \{ z_i=(x_i, u_i), x^+_i, r_i \}_{i \in [1:n]}$; 
	 discount factor~$\gamma \in (0,1)$; 
	 MP features~$f:\Z \ra \Rl^p$; 
	 termination constant~$\epsilon$.
	\ENSURE MP coefficients~$\theta \in \R^p$.

  	\STATE compute $F\sam,\ G\sam \in  \Rl^{n\times p}$ as in Proposition~\ref{prop:MP_Bellman};
  	\STATE initialize by $\theta \gets (0,\ldots,0)\tr \in \R^p$; 
   $\theta^{+} \gets \mathbf{D}\theta$ as in \eqref{eq:MP-FQI_iteration};
  	
  	\WHILE{$\norm{\theta^+ - \theta}_{\infty} > \epsilon$}
  	
  	\STATE $\theta \gets \theta^+$; 
   $\theta^{+} \gets \mathbf{D}\theta$ as in \eqref{eq:MP-FQI_iteration};
  	
  	\ENDWHILE 
  	\STATE output $\theta \gets \theta^+$.
\end{algorithmic}
\end{small}
\end{algorithm}

\subsection{Analysis}
Notice that Algorithm~\ref{alg:mp-fqi} terminates when the difference between the output of two consecutive iterations (in $\infty$-norm) is less than a prescribed tolerance $\epsilon > 0$. 
As shown by the following result, such a condition is sufficient for Algorithm~\ref{alg:mp-fqi} to terminate in a finite number of iterations. 

\begin{Thm} [Convergence of MP-FQI] \label{thm:mp-fqi conv}
Let Assumption~\ref{As:features} hold. Then, the operator $\mathbf{D}:\R^p \ra \R^p$ defined in~\eqref{eq:MP-FQI_iteration} is a $\gamma$-contraction in $\infty$-norm. 
\end{Thm}

The preceding result shows that Algorithm~\ref{alg:mp-fqi} converges linearly with a rate $\leq \gamma$. 
Moreover, in combination with the result of Proposition~\ref{prop:MP_FQI_reg}, it implies that the recursion~\eqref{eq:MP-FQI_iter0} converges to $\theta^{(\infty)} = \tilde{\theta}^{(\infty)} + \beta^{(\infty)} e$, where $\tilde{\theta}^{(\infty)} = \mathbf{D} \tilde{\theta}^{(\infty)}$ is the unique fixed-point of $\mathbf{D}$ and 
$$\beta^{(\infty)} = \frac{1}{2(1-\gamma)} \Vert{F\sam} \boxtimes \tilde{\theta}^{(\infty)} -  G\sam \boxtimes (\gamma\tilde{\theta}^{(\infty)})\Vert_{\infty}.$$
We next consider the complexity of the proposed algorithm.

\begin{Thm} [Complexity of MP-FQI] \label{thm:mp-fqi comp}
For Algorithm~\ref{alg:mp-fqi}, disregarding the complexity of the maximization over $u^+$ for computing $G\sam$, the compilation and per-iteration time complexities are both of $\ord(np)$. 
\end{Thm}

Let us also discuss the complexity of solving the maximization over the action~$u^+$ in Algorithm~\ref{alg:mp-fqi}. 
Note that these maximizations are to be solved \emph{only once}. 
In particular, if the action space $\U = \{v_1, \ldots, v_{p_{\mathrm{u}}} \}$ is finite and
\begin{equation}\label{eq:feature_sep}
   [f]_{jk}(x,u) = [f_{\mathrm{x}}]_j(x)+\delta_{v_k} (u), \quad j \in [1:p_{\mathrm{x}}],\ k \in [1:p_{\mathrm{u}}], 
\end{equation}
with state features $f_{\mathrm{x}}:\X\ra\Rl^{p_{\mathrm{x}}}$, we have
\begin{align*}
 \max_{u^+} [f]_{jk}(x^+_i, u^+) &=  [f_{\mathrm{x}}]_j(x^+_i).
\end{align*}
That is, we do \emph{not} need to solve any maximization problem.

We finish this section by comparing the proposed MP-FQI algorithm with MP-linear parametrization with its counterpart in the conventional plus-times algebra. 

\begin{Rem}[Comparison with standard FQI] \label{rem:fqi_trad}
Consider the FQI algorithm given by~\eqref{eq:fqi_stand} with \emph{linear} parametrization $Q_{\theta}(z) = \Tilde{f}(z)\tr \theta$, where $\Tilde{f}:\Z \ra \R^p$ is the vector of features, and \emph{ridge} regression~\cite{Hoerl}. 
For these choices of parametrization and loss, the recursion~\eqref{eq:fqi_stand} reduces to
\begin{equation}\label{eq:fqi LS}
\theta^{(\ell+1)}  \in \argmin_{\theta}\ J(\theta) + \lambda \norm{\theta}_2^2, \quad \ell\in\N_0,
\end{equation}  
where
\[
J(\theta) = \sum_{i\in [1:n]} \bigg( \Tilde{f}(z_i)\tr \theta  - \big( r_i + \gamma \cdot \max_{u^+} \Tilde{f}(x^+_i, u^+)\tr \theta^{(\ell)}  \big)\bigg)^2,
\]
and $\lambda\geq 0$ is the regularization coefficient. 
This is a convex, quadratic program. 
By defining the matrix $\Tilde{F}\sam\in \R^{n\times p}$ with entries $[\Tilde{F}\sam]_{i,j} = [\Tilde{f}]_j(z_i)$, the recursion~\eqref{eq:fqi LS} can be written as 
\begin{equation}\label{eq:fqi_stand_1}
\left\{\begin{array}{l}
[g\sam^{(\ell)}]_i = r_i + \gamma \cdot \max_{u^+} \Tilde{f}(x^+_i, u^+)\tr \theta^{(\ell)},\ i \in [1:n], \\
\theta^{(\ell+1)}  = \big(\Tilde{F}\sam\tr \Tilde{F}\sam + \lambda I\big)^{-1}\Tilde{F}\sam\tr g\sam^{(\ell)}.
\end{array}\right.
\end{equation}
We note that (i)~this algorithm is \emph{not} guaranteed to converge, and (ii)~the time complexity of compilation is of $\ord(np^2+p^3)$ while each iteration requires $\ord(np)$ operations besides solving a maximization over $u^+$ for each $i \in [1:n]$.   
\end{Rem}

\subsection{An alternative implementation of MP-FQI}
We now present an alternative characterization of the solution to the empirical BE~\eqref{eq:DP_empirical_2} 
which is inspired by the MP eigenvector method developed for $H_{\infty}$ control problems~\cite{McEneaney03, McEneaney04}. 
As we see shortly, this characterization leads to an algorithm with a per-iteration complexity independent of the sample size~$n$ at the cost of increasing the compilation time. 

\begin{Prop}[MP empirical BE II]\label{prop:MP_Bellman II}
Consider the matrices $F\sam$ and $G\sam$ in Proposition~\ref{prop:MP_Bellman}. 
Let
\begin{itemize}
    \item the matrix $\bar{C}\sam \in \Rl^{p\times p}$ be a solution to the MP-linear equation $F\sam \boxtimes C\sam = G\sam$, and, 
    \item the vector $\bar{\theta} \in \R^p$ be a solution to the MP-fixed-point equation $\theta = \bar{C}\sam \boxtimes (\gamma \theta)$.
\end{itemize}
Then, the MP-linear function $Q_{\bar{\theta}}(z) = f(z)\tr \boxtimes \bar{\theta}$ solves the empirical BE~\eqref{eq:DP_empirical_2}. 
\end{Prop}

Let us also note that the equality $F\sam \boxtimes C\sam = G\sam$ implies that 
\begin{equation*}\label{eq:Theta def app}
\sdp [f]_j (z_i) = \max_{k \in [1:p]} \left\{ [f]_k(z_i) + [C\sam]_{k,j} \right\}, \quad i \in [1:n],\ j \in [1:p], 
\end{equation*}
that is, the $j$-th column of $C\sam$ contains the max-plus coefficients of the MP-linear representation of $\sdp [f]_j$ \emph{with respect to the same features $f:\Z\ra\Rl^p$}, at the sample points in $\Z\sam$. 

The preceding lemma points to a recursive algorithm based on the fixed-point iteration 
\begin{equation}\label{eq:fp_MPFQI}
   \theta^{(\ell +1)} = C\sam \boxtimes (\gamma \theta^{(\ell)}),\quad \ell \in \N_0, 
\end{equation}
with some initialization $\theta^{(0)} \in \R^p$.
What remains to be addressed is the computation of the matrix $C\sam \in \Rl^{p\times p}$, i.e., a solution to the MP-linear equation $F\sam \boxtimes C\sam = G\sam$. 
A possible choice is the solution to the column-wise regression problems
\begin{align*}
   [\wh{C}\sam]_{\cdot,j} &= \argmin_{c} \norm{F\sam \boxtimes c - [G\sam]_{\cdot,j}}_{\infty} 
   = c_{\text{ps}}^{(j)} + \frac{\norm{F\sam\boxtimes c_{\text{ps}} - [G\sam]_{\cdot,j}}_{\infty}}{2}  e, 
\end{align*}
where 
\[
c_{\text{ps}}^{(j)} = - \left( F\sam\tr \boxtimes (-[G\sam]_{\cdot,j}) \right),
\]
for each $j \in [1:p]$. 
This choice of the matrix $C\sam$ leads to a \emph{compilation} time complexity of $\ord(np^2)$. 
On the other hand, 
the \emph{per-iteration} time complexity in~\eqref{eq:fp_MPFQI} is of $\ord(p^2)$ which is independent of the sample size. 
We also note that the iteration~\eqref{eq:fp_MPFQI} converges linearly with rate $\gamma$ assuming that the matrix~$C\sam$ is row $\R$-astic; see Lemma~\ref{lem:contraction general}.

\section{Variational max-plus FQI}
\label{sec:variational}

In this section, we propose and analyze the \emph{variational} implementation of the max-plus fitted Q-iteration (v-MP-FQI). 
The proposed algorithm is based on an alternative formulation of the BE which originates from the MP variational formulation of the continuous-time, deterministic optimal control problems proposed in~\cite{Akian08}. 
This formulation is also used for approximating the value function in the optimal control problem of deterministic MDPs by~\cite{Bach20}. 
To this end, let us define the MP ``inner product'' of the two functions $f,\tilde{f} \in \Rl^{\Z}$ over the state-action space $\Z$ by 
\[
f \boxdot \tilde{f} \Let \max_{z \in \Z} \{ f(z) + \tilde{f}(z) \}.
\]
Then, for the vector of \emph{test functions}~$h: \Z \ra \Rl^q $, we consider finding $Q\opt$ by solving the \emph{variational} BE
\begin{equation}\label{eq:var Bellman}
[h]_k \boxdot Q = [h]_k \boxdot \dpo Q, \quad \forall k \in [1:q].
\end{equation}

\subsection{Algorithm} 
Similar to the MP-FQI algorithm, we incorporate the MP-linear approximation~\eqref{eq:mp-approx Q} for the Q-function. 
However, instead of the empirical BE~\eqref{eq:DP_empirical_2}, we aim to solve the empirical version of the \emph{variational} BE~\eqref{eq:var Bellman}, given by 
\begin{equation}\label{eq:var Bellman sampled}
[h]_k \boxdot\sam Q_{\theta}  = [h]_k \boxdot\sam \sdp Q_{\theta} , \quad \forall k \in [1:q],
\end{equation}
where  
\[
f \boxdot\sam \tilde{f} \Let \max_{i\in[1:n]} \{ f(z_i) + \tilde{f}(z_i) \},
\]
is the \emph{empirical} MP inner product of $f,\tilde{f} \in \Rl^{\Z}$ computed over the sampled state-action space $\Z\sam$. 
Once again, we provide an alternative characterization of the preceding equation which forms the basis of the proposed algorithm.

\begin{Prop}[MP empirical variational BE]\label{prop:MP_Bellman_var}
Consider the matrices $F\sam$ and $G\sam$ in Proposition~\ref{prop:MP_Bellman}, and 
define the matrix $H\sam \in \Rl^{n\times q}$ with entries
\begin{equation*}
    [H\sam]_{i,k} = [h]_k(z_i), \qquad i\in[1:n],\ k\in[1:q]. 
\end{equation*}
Set $F\sam^{H} = H\sam\tr  \boxtimes F\sam  \in \Rl^{q\times p} $ and $G\sam^{H} = H\sam\tr \boxtimes G\sam  \in \Rl^{q\times p} $, and let $\bar{\theta} \in \R^p$ be a solution to the MP equation
\begin{equation}\label{eq:DP_empirical_mp_var}
    F\sam^{H} \boxtimes \theta = G\sam^{H} \boxtimes (\gamma\theta).
\end{equation}
Then, the MP-linear function $Q_{\bar{\theta}}(z) = f(z)\tr \boxtimes \bar{\theta}$ in \eqref{eq:mp-approx Q} solves the empirical variational BE~\eqref{eq:var Bellman sampled}. 
\end{Prop}

Using the preceding result, once again we see that the parameters of the MP-linear approximation~\eqref{eq:mp-approx Q} of the Q-function can be determined by solving the MP equation~\eqref{eq:DP_empirical_mp_var}. 
We next follow a similar procedure to that of the previous section for solving \eqref{eq:DP_empirical_mp_var}. 
To this end, we need the following assumption on the samples and the test functions:

\begin{As}
\label{As:compatible test} 
The test functions~$h:\Z\ra\Rl^q$ and the samples~$ \{ z_i=(x_i, u_i), x^+_i, r_i \}_{i \in [1:n]}$ are such that 
\begin{align}
    &\forall z \in \Z,\ \exists k \in [1:q],\ [h]_k(z) \neq -\infty; \label{eq:test proper}  \\
    &\forall k \in [1:q],\ \exists i \in [1:n],\ [h]_k(z_i) \neq -\infty. \label{eq:H col realistic}
\end{align}
That is, $H\sam$ is doubly $\R$-astic.
\end{As}

The preceding assumption is relevant if the tests functions are \emph{extended} real-valued functions and assign $-\infty$ to a subset of the state-action space~$\Z = \X\times\U$. 
The implications of this assumption are also similar to that of Assumption~\ref{As:features}: 
It requires the test functions to be chosen properly and the samples to be rich enough. 
To be precise, condition~\eqref{eq:test proper} implies that each state-action pair $z = (x,u)$ activates at least one test function. 
If this is not the case, then the test functions are not chosen properly in the sense that any max-plus linear combination of them assigns $-\infty$ to a subset of $\Z$. Condition~\eqref{eq:H col realistic}, on the other hand, implies that each test function $[h]_k$ is activated by at least one sample $z_i$. 
If this condition does not hold for some $k\in[1:q]$, the corresponding test function $[h]_k$ can be removed since both sides of the empirical variational BE~\eqref{eq:var Bellman sampled} are equal to $-\infty$. 

Assumption~\ref{As:compatible test} combined with Assumption~\ref{As:features}, leads to $F\sam^{H}$ and $G\sam^{H}$ having the following desired properties:

\begin{Lem} \label{lem:mp-vqi conv}
Under Assumptions~\ref{As:features} and \ref{As:compatible test}, $F\sam^{H}$ is doubly $\R$-astic and $G\sam^{H}$ is row $\R$-astic.
\end{Lem}

Therefore, under Assumptions~\ref{As:features} and \ref{As:compatible test}, we can follow a similar procedure to that of the previous section for the MP-FQI algorithm and use recursive unconstrained MP-linear regressions in the $\infty$-norm for solving \eqref{eq:DP_empirical_mp_var} while only keeping track of the iterations of the principal solution; cf.~\eqref{eq:MP-FQI_iter0} and \eqref{eq:MP-FQI_iteration}. 
This leads to the recursion
\begin{equation}\label{eq:MP-FQI_iteration_var}
\theta^{(\ell+1)} = \mathbf{D}^{H} \theta^{(\ell)} \Let  - \left[ (F^{H}\sam)\tr \boxtimes \big[-(G^{H}\sam \boxtimes (\gamma\theta^{(\ell)}))\big] \right],  
\end{equation}
for the MP coefficients~$\theta$.  
Algorithm~\ref{alg:mp-fqi-var} summarizes the procedure described above.

\begin{algorithm}[t]
\begin{small}
   \caption{Variational max-plus FQI (v-MP-FQI)} 
   \label{alg:mp-fqi-var}
\begin{algorithmic}[1]
	\REQUIRE samples~$ \{ z_i=(x_i, u_i), x^+_i, r_i \}_{i \in [1:n]}$; 
	 discount factor~$\gamma \in (0,1)$; 
	 MP features~$f:\Z \ra \Rl^p$;
      MP test functions~$h:\Z \ra \Rl^q$;
	 termination constant~$\epsilon$.
	\ENSURE MP coefficients~$\theta \in \R^p$.

   \STATE compute $F\sam^{H},\ G\sam^{H} \in \Rl^{q\times p} $ as in Propositions~\ref{prop:MP_Bellman} and \ref{prop:MP_Bellman_var};
  	\STATE initialize by $\theta \gets (0,\ldots,0)\tr \in \R^p$; $\theta^{+} \gets \mathbf{D}^{H} \theta $ as in~\eqref{eq:MP-FQI_iteration_var}; 
  	
  	\WHILE{$\norm{\theta^+ - \theta}_{\infty} \ge \epsilon$}
  	
  	\STATE $\theta \gets \theta^+$; 
   $\theta^{+} \gets \mathbf{D}^{H} \theta $ as in~\eqref{eq:MP-FQI_iteration_var};
  	
  	\ENDWHILE 
  	\STATE output $\theta \gets \theta^+$.
\end{algorithmic}
\end{small}
\end{algorithm}

\subsection{Analysis} 
We now consider the convergence and complexity of Algorithm~\ref{alg:mp-fqi-var}. 

\begin{Thm} [Convergence of v-MP-FQI] \label{thm:mp-vqi conv}
Let Assumptions~\ref{As:features} and \ref{As:compatible test} hold. 
Then, the operator~$\mathbf{D}^{H}:\R^p \ra \R^p$ defined in~\eqref{eq:MP-FQI_iteration_var} is a $\gamma$-contraction in $\infty$-norm.  
\end{Thm}

\begin{Thm} [Complexity of v-MPFQI] \label{thm:mp-vqi comp}
For Algorithm~\ref{alg:mp-fqi-var}, disregarding the complexity of the maximization over $u^+$ for computing $G\sam$, the compilation time complexity is of $\ord(npq)$. 
The per-iteration complexity of Algorithm~\ref{alg:mp-fqi-var} is of $\ord(pq)$. 
\end{Thm}

Therefore, the v-MP-FQI Algorithm~\ref{alg:mp-fqi-var} also converges linearly with a rate $\leq \gamma$. However, notice how the \emph{per-iteration} time complexity of Algorithm~\ref{alg:mp-fqi-var} is independent of the sample size and instead depends on the number of test functions; cf.~Theorem~\ref{thm:mp-fqi comp} for the complexity of Algorithm~\ref{alg:mp-fqi}. 
The important property of v-MP-FQI is hence that it goes through the samples only once. 
Hence, for a large sample size~$n$ (which requires a lot of computational resources per iteration of the (MP-)FQI algorithm), 
the v-MP-FQI algorithm can be employed to control the time requirement by choosing a proper number~$q$ of test functions. We finish this section with a remark on the computation of $F\sam^H$.

\begin{Rem}[Exact computation of $F^H$] 
We might be able to compute the MP inner product on the right-hand-side of the variational BE~\eqref{eq:var Bellman} exactly. 
Given features $f:\Z \ra \Rl^p$ and test functions $h:\Z \ra \Rl^q$, we have 
$$ \max_{z \in \Z} \{ [h]_k(z) + Q_{\theta}(z) \} = [F^{H} \boxtimes \theta]_k, \quad k \in [1:q],$$
where $F^{H}  \in \Rl^{q\times p}$ is the exact version of $F\sam^H$ with entries 
\begin{align*}
[F^{H}]_{k,j} =  \max_{z \in \Z} \left\{ [h]_k(z) + [f]_j(z)  \right\}.
\end{align*} 
If $F^{H}$ is available, we can use it instead of $F\sam^H$ in Algorithm~\ref{alg:mp-fqi-var}. 
For instance, for quadratic features $[f]_j(z) = -\frac{c}{2}\norm{z-w_j}_2^2$ and test functions $[h]_k(z) = -\frac{\tilde{c}}{2}\norm{z-\tilde{w}_k}_2^2$, with $\{w_j\}_{j \in [1:p]},\{\tilde{w}_k\}_{k \in [1:q]} \subset \Z$ and $c,\tilde{c}>0$, we have $[F^{H}]_{k,j} = \frac{cw_j+\tilde{c}\tilde{w}_k}{c+\tilde{c}}$, assuming $\Z$ is convex. 
\end{Rem}

\section{Numerical experiments}
\label{sec:num ex}
In this section, we examine the performance of the proposed MP-FQI Algortihm~\ref{alg:mp-fqi} and v-MP-FQI Algorithm~\ref{alg:mp-fqi-var} in comparison with the (standard) FQI algorithm~\eqref{eq:fqi_stand_1}. 
To this end, we consider the DC motor stabilization problem adapted from~\cite[Sec.~3.4.5]{Busoniu10}. 
A detailed description of the problem and setup of the numerical simulations is provided in Appendix~\ref{app:num}. 
In particular, we consider a discrete action space~$\U= \{ v_k = -10+5(k-1)\}_{k \in [1:5]}$ which is the discretization of the continuous action space~$[-10,10]$. 
The sample state-action pairs~$\{z_i = (x_i,u_i)\}_{i \in [1:n]}$, with $n=5000$, are generated uniformly at random from the state-action space~$\Z = \X\times \U$. 
For MP-FQI and v-MP-FQI, we combined (i)~quadratic and (ii)~distance \emph{state} features with MP-indicator \emph{input} features~$\delta_{v_k}$ as in~\eqref{eq:feature_sep} with the total number of features being $p = p_{\mathrm{u}} \cdot p_{\mathrm{x}} = 5p_{\mathrm{x}}$. 
The test functions in v-MP-FQI are chosen to be the same as the features so that $q = p$. 
For FQI, we similarly combine (i)~radial basis functions (RBFs) and (ii)~indicator functions as the \emph{state} features with indicator \emph{input} features~$\ind{v_k}$ with the same total number of features $p = p_{\mathrm{u}} \cdot p_{\mathrm{x}} = 5p_{\mathrm{x}}$. 
The regularization parameter in \eqref{eq:fqi_stand_1} is set to $\lambda = 10^{-3}$.
Finally, for termination of the algorithms, besides a maximum number of iterations ($\ell_{\max} = 1000$), we use the condition $\norm{\theta^+ - \theta}_{\infty} \leq 10^{-3} \norm{\theta}_{\infty}$ in order to create a fair comparison across algorithms that generates $\theta$ in different orders of magnitude.

\begin{figure*}
\centering
\includegraphics[width=.49\linewidth]{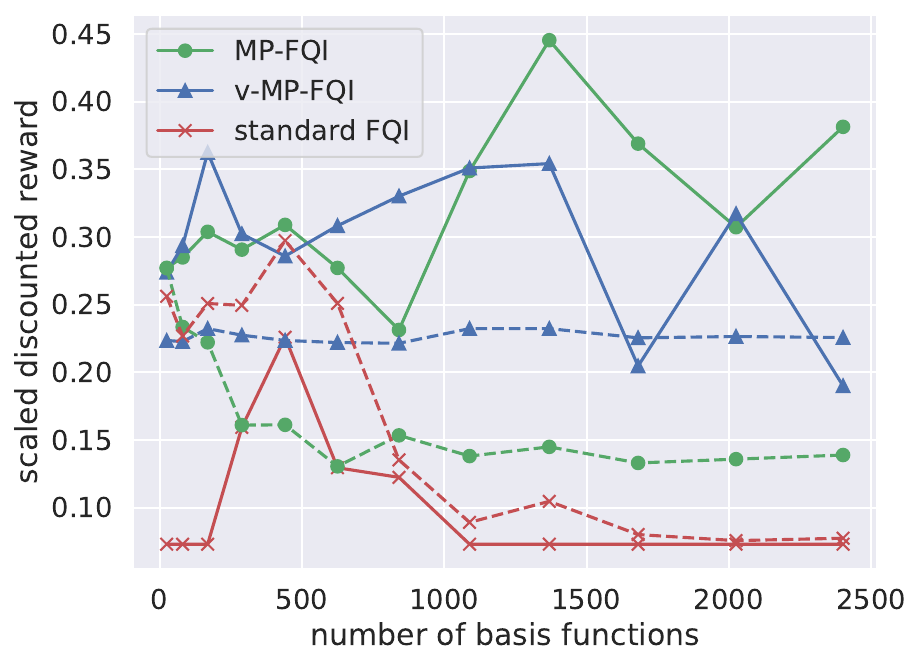}
\includegraphics[width=.49\linewidth]{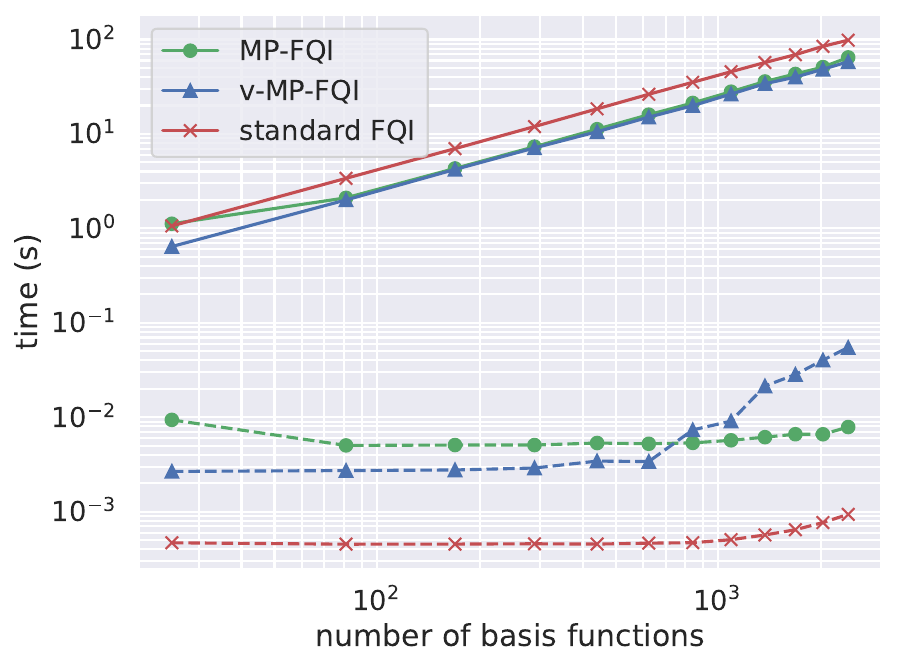}
\includegraphics[width=.49\linewidth]{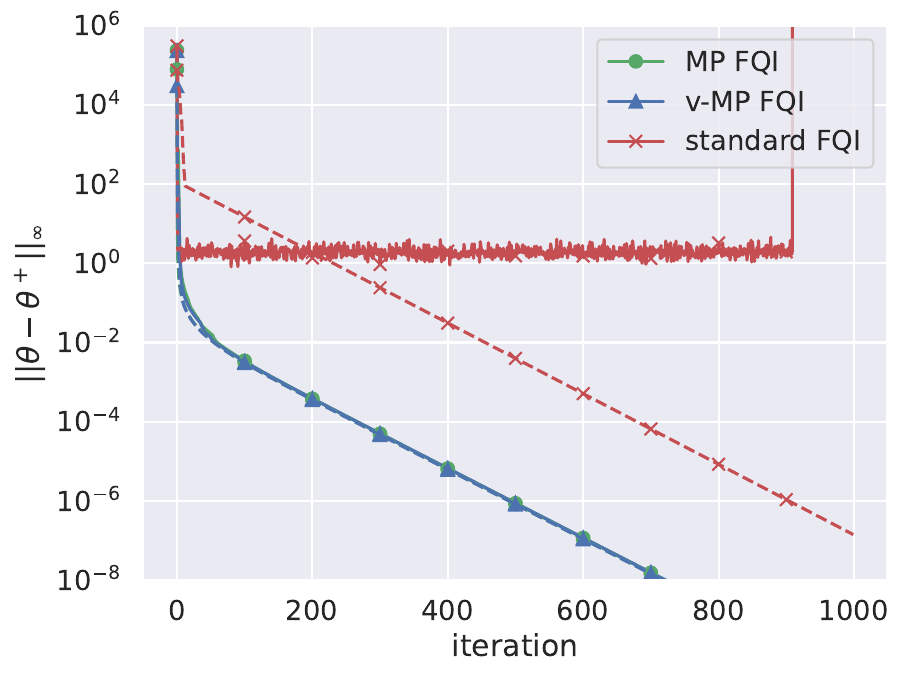}
\caption{DC motor stabilization problem. \emph{Top-Left}:~Average reward of 100 instances of the problem with random initial state over $T=100$ time steps. Solid (resp.~dashed) lines correspond to quadratic (resp.~distance) functions in (v-)MP-FQI and RBF (resp.~indicator functions) in FQI for state features. \emph{Top-Right}:~The running time of the algorithms. Solid lines (resp.~dashed) correspond to compilation (resp.~per-iteration) time. \emph{Bottom}:~Convergence of algorithms. Solid and dashed lines are the same as in the top-left figure.}
\label{fig:DC general}
\end{figure*}

The left panel of Figure~\ref{fig:DC general} shows the average reward of 100 instances of the problem for the greedy control policy~$\pi(x) \in \argmin_{u \in \U }  Q_{\theta} (x,u)$, 
using the output~$\theta$ of these algorithms. 
The reported rewards are normalized against the reward of the linear quadratic controller. 
As can be seen, the proposed algorithms, MP-FQI and v-MP-FQI, consistently outperform their standard counterpart when it comes to the quality of the greedy policy. 

The compilation and per-iteration running time of the algorithms is reported in the middle panel of Figure~\ref{fig:DC general}. 
Some remarks are in order regarding the reported running times: 
The dominant factor in the compilation time is the computation of the matrices $F\sam$, $G\sam$, and $H\sam$. 
In particular, since we chose the same test functions as the features in v-MP-FQI and hence $H\sam = F\sam$, the compilation times of MP-FQI and v-MPFQI are almost the same despite the extra MP matrix multiplication required for computing $F\sam^{H}$ and $G\sam^{H}$ in v-MPFQI. 
We note that the run-times reported here are derived from single-thread operations, however, the computation of these matrices is embarrassingly parallelizable. 
On the other hand, the reported per-iteration run-times match the complexities given in Theorems~\ref{thm:mp-fqi comp} and \ref{thm:mp-vqi comp}. 
In particular, when $p=q <n$, v-MP-FQI has a lower time requirement compared to MP-FQI. 
However, as $p=q$ increases and becomes larger than $n$, the run-time of v-MP-FQI increases with a quadratic rate and surpasses that of MP-FQI. 

Finally, in the right panel of Figure~\ref{fig:DC general}, we report $\norm{\theta^+ - \theta}_{\infty}$ in successive iterations to examine the convergence of the algorithms. 
As shown, corresponding to Theorems~\ref{thm:mp-fqi conv} and~\ref{thm:mp-vqi conv}, the proposed MP-FQI and v-MP-FQI are convergent. 
However, the convergence of the standard FQI with a linear approximation is not in general guaranteed. As we can see, FQI diverges with the global RBFs but converges with the local indicator functions as the state features.

\section{Final remarks}
\label{sec:conclusion}

In this paper, we introduced a novel FQI algorithm based on an MP-linear approximation of the Q-function with a provable convergence. 
We also considered the variational implementation of the algorithm leading to a per-iteration complexity independent of the number of samples. 
We now provide some final remarks on the limitations and extensions of the proposed algorithms and the provided analysis.  

One of the main drawbacks of the current work is the lack of error analysis, i.e., a theoretical bound on the distance between the output of the algorithms and the optimal Q-function. 
This issue can be addressed by bounding the error introduced in solving the MP-linear regression problem with respect to the output of the true Bellman operator in each iteration~\cite{bertsekas1996neuro}. 
Also of interest is the identification of the particular class of problems (e.g., regularity assumptions on the transition kernel and reward function of the MDP) for which the optimal value function can be well-represented in MP-linear function spaces. 
A classic example is the case of linear dynamics with additive disturbance with a concave reward function which leads to a concave optimal Q-function~\cite{bert1973lin}. 

A possible extension of the current algorithms is their potential combination with fast numerical algorithms for discrete conjugate or distance transforms to reduce the computation cost of the algorithm. 
In particular, with factorized (i.e., grid-like) samples and/or features, we can use a similar approach used in~\cite{Lucet09, kolari2023TAC} to reduce the per-iteration complexity of Algorithm~\ref{alg:mp-fqi} to $\ord(n+p)$. 
Another interesting extension is the application of sparse approximate solutions to MP-linear equations~\cite{tsi2021spa} instead of the MP-linear regression in the proposed algorithms. 
Such sparse solutions allow the user to control the size of the approximation (i.e., number of parameters) and reduce the memory and time complexity required for generating the greedy policy.

\appendix

\section{Technical proofs}\label{app:proof}
\subsection{Proof of Lemma~\ref{lem:contraction general}}
We first note that the assumption on matrix $A$ being row $\R$-astic implies that the range of $\mathbf{A}$ is indeed (a subset of) $\R^m$. 
Now, let $y_1, y_2 \in \R^n$ and observe that for each $i \in [1:m]$, we have
\begin{align*}
[\mathbf{A}y_2]_i &=[A \boxtimes (\eta y_2)]_i = \max_{j \in [1:n]} \left\{ [A]_{i,j} + \eta \cdot [y_2]_j \right\} 
= \max_{j \in [1:n]} \left\{ [A]_{i,j} + \eta \cdot [y_1]_j + \eta \left( [y_2]_j - [y_1]_j\right) \right\}.
\end{align*}
Hence, using the fact that 
$$[\mathbf{A}y_1]_i =  \max_{j \in [1:n]} \left\{ [A]_{i,j} + \eta \cdot [y_1]_j \right\},$$ 
we have 
\begin{align*}
[\mathbf{A}y_1]_i - \eta  \norm{y_1 - y_2}_{\infty} \leq [\mathbf{A}y_2]_i \leq [\mathbf{A}y_1]_i + \eta  \norm{y_1 - y_2}_{\infty}.
\end{align*}
That is, 
\begin{align*}
\big|[\mathbf{A}y_1]_i - [\mathbf{A}y_2]_i \big| \leq \eta  \norm{y_1 - y_2}_{\infty} , \quad \forall i \in [1:m]. 
\end{align*}
This completes the proof. 
\subsection{Proof of Lemma~\ref{lem:mp add and hom}}

Let $z_i \in \Z\sam$. 
We have
\begin{align*}
\sdp [\max \{ f, \tilde{f} \}](z_i) 
&= \rew_i + \gamma \cdot \max_{u^+} \ [\max \{ f, \tilde{f} \}] (x_i^+,u^+) \\
&= \rew_i + \gamma \cdot \max_{u^+} \ \max \left\{ f(x_i^+,u^+),\tilde{f}(x_i^+,u^+) \right\} \\
&= \rew_i + \gamma \cdot \max \big\{  \max_{u^+}  f(x_i^+,u^+), \max_{u^+} \tilde{f}(x_i^+,u^+) \big\} \\
&= \max \big\{ \rew_i + \gamma \cdot  \max_{u^+}  f(x_i^+,u^+), \rew_i + \gamma \cdot \max_{u^+} \tilde{f}(x_i^+,u^+) \big\} \\
&= \max \big\{ \sdp f (z_i), \sdp \tilde{f}(z_i) \big\} \\
&= [\max \{ \sdp f , \sdp \tilde{f} \}] (z_i).
\end{align*}
Similarly, 
\begin{align*}
\sdp [\alpha + f](z_i) &= \rew_i + \gamma \cdot \max_{u^+} \ [\alpha + f] (x_i^+,u^+) \\
&= \rew_i + \gamma \cdot \max_{u^+} \left\{ f(x_i^+,u^+) + \alpha \right\} \\
&= \gamma  \alpha +  \rew_i +  \gamma \cdot \max_{u^+}  f(x_i^+,u^+) \\
&= \gamma  \alpha  +  \sdp f (z_i) = [\gamma  \alpha  + \sdp f] (z_i).
\end{align*}
\subsection{Proof of Proposition~\ref{prop:MP_Bellman}} 
Fix $i\in [1:n]$. From the definition of $F\sam$, we have
\begin{align*}
    Q_{\bar{\theta}}(z_i) & =  f(z_i)\tr \boxtimes \bar{\theta} 
    =  \max_{j\in [1:p]} \left\{ [f]_j(z_i) + [\bar{\theta}]_j \right\}  
    =  \max_{j\in [1:p]} \left\{ [F\sam]_{i,j} + [\bar{\theta}]_j \right\} 
    =  [{F\sam} \boxtimes \bar{\theta}]_i.
\end{align*}
Hence, using \eqref{eq:DP_empirical_mp} and the definition of $G\sam$, we can write
\begin{align*}
    Q_{\bar{\theta}}(z_i)  &= [G\sam \boxtimes (\gamma\bar{\theta})]_i 
     =  \max_{j\in [1:p]} \left\{ [G\sam]_{i,j} + \gamma  [ \bar{\theta}]_j \right\}  
     =  \max_{j\in [1:p]} \left\{ \sdp [f]_j(z_i) + \gamma [ \bar{\theta}]_j \right\}. 
\end{align*}
Finally, we can use Lemma~\ref{lem:mp add and hom} to obtain
\begin{align*}
    Q_{\bar{\theta}}(z_i)  &=  \max_{j\in [1:p]} \left\{ \sdp \big[ [f]_j+ [\bar{\theta}]_j \big](z_i) \right\} 
    =   \sdp \big[ \max_{j\in [1:p]} \left\{ [f]_j + [\bar{\theta}]_j \right\} \big](z_i)  
     =   \sdp Q_{\bar{\theta}}(z_i).
\end{align*}
This completes the proof.

\subsection{Proof of Proposition~\ref{prop:MP_FQI_reg}}
The proof is by induction. The base case $\ell = 0$ holds by construction. 
Assume now $\theta^{(\ell)} = \tilde{\theta}^{(\ell)} + \beta^{(\ell)} e$. 
Then, under Assumption~\ref{As:features}, the solution to the MP-linear regression~\eqref{eq:MP-FQI_iter0} is given by~\cite[Thm.~3.5.2]{butkovivc10}
\begin{equation*}
    \theta^{(\ell+1)} = y + \alpha e  , \quad \ell \in \N_0,
\end{equation*}
where
\begin{equation}\label{eq:reg_ps}
y =   - \left[ F\sam\tr \boxtimes \big[-(G\sam \boxtimes (\gamma \theta^{(\ell)}))\big] \right],  
\end{equation}
is the \emph{principal solution} of ${F\sam} \boxtimes \theta =  G\sam \boxtimes (\gamma\theta^{(\ell)})$ and 
\begin{equation}\label{eq:reg_const}
\alpha = \frac{1}{2} \Vert{F\sam} \boxtimes y -  G\sam \boxtimes (\gamma\theta^{(\ell)})\Vert_{\infty},
\end{equation}
is the constant shift in the solution. 
Now, observe that $M\boxtimes(v + c e ) = M\boxtimes v + c e $ for matrix~$M$, vector~$v$ and scalar~$c$. 
Then, by plugging $\theta^{(\ell)} = \tilde{\theta}^{(\ell)} + \beta^{(\ell)} e$ in \eqref{eq:reg_ps}, we have
\begin{align*}
   y & = - \left[ F\sam\tr \boxtimes \big[-\big(G\sam \boxtimes (\gamma (\tilde{\theta}^{(\ell)} + \beta^{(\ell)} e))\big)\big] \right] \\
   & = - \left[ F\sam\tr \boxtimes \big[-\big(G\sam \boxtimes (\gamma \tilde{\theta}^{(\ell)})\big)\big] \right] + \gamma \beta^{(\ell)} e\\
   & =   \mathbf{D}\tilde{\theta}^{(\ell)} + \gamma \beta^{(\ell)} e = \tilde{\theta}^{(\ell+1)} + \gamma \beta^{(\ell)} e. 
\end{align*}
Similarly, for~\eqref{eq:reg_const}, we have
\begin{align*}
    \alpha &= \frac{1}{2} \Vert{F\sam} \boxtimes (\tilde{\theta}^{(\ell+1)} + \gamma \beta^{(\ell)} e) -  G\sam \boxtimes \big(\gamma (\tilde{\theta}^{(\ell)} + \beta^{(\ell)} e)\big)\Vert_{\infty} \\
    & = \frac{1}{2} \Vert{F\sam} \boxtimes \tilde{\theta}^{(\ell+1)} -  G\sam \boxtimes (\gamma\tilde{\theta}^{(\ell)})\Vert_{\infty}.
\end{align*}
Therefore, we have 
$$\theta^{(\ell+1)} = \tilde{\theta}^{(\ell+1)} + \gamma \beta^{(\ell)} e +\alpha e = \tilde{\theta}^{(\ell+1)} + \beta^{(\ell+1)} e.$$ 
\subsection{Proof of Theorem~\ref{thm:mp-fqi conv}}

Let $\theta_1,\theta_2 \in \R^p$ and 
recall that Assumption~\ref{As:features} imply that the matrices $F\sam\tr$ and $G\sam$ are row $\R$-astic. 
This, in turn, implies that $(G\sam \boxtimes \gamma\theta_1) \in \R^p$ and $(G\sam \boxtimes \gamma\theta_2) \in \R^p$. 
Hence, we can use Lemma~\ref{lem:contraction general} with $A=F\sam\tr$ and $\eta = 1$ to write
\begin{align*}
\norm{\mathbf{D}\theta_1 - \mathbf{D}\theta_2}_{\infty} 
&= \norm{F\sam\tr \boxtimes \left(-(G\sam \boxtimes (\gamma\theta_1))\right) -  F\sam\tr \boxtimes \left(-(G\sam \boxtimes (\gamma\theta_2))\right) }_{\infty} \\
&\leq \norm{ G\sam \boxtimes (\gamma\theta_1) - G\sam \boxtimes (\gamma\theta_2)}_{\infty}.
\end{align*} 
Another application of Lemma~\ref{lem:contraction general} with $A=G\sam$ and $\eta = \gamma$ then gives us $\norm{\mathbf{D}\theta_1 - \mathbf{D}\theta_2}_{\infty} \leq \gamma  \norm{\theta_1 - \theta_2}_{\infty}$.
\subsection{Proof of Theorem~\ref{thm:mp-fqi comp}}

The compilation time includes the computation of the matrices $F\sam$ and $G\sam$ with $np$ entries, which requires $\ord(np)$ operations. 
Within each iteration, the parameter update involves a scalar multiplication of size $n$ and two MP matrix-vector multiplications of size $n \times p$ and $p \times n$. 
These computations also require $\ord(np)$ operations. 
\subsection{Proof of Proposition~\ref{prop:MP_Bellman II}}

Fix $i \in [1:n]$. Observe that 
\begin{align*}
   \sdp Q_{\bar{\theta}}(z_i) = \sdp \left[  f\tr \boxtimes \bar{\theta} \right](z_i)
   =  \sdp \big[ \max_{j\in [1:p]} \left\{ [f]_j + [\bar{\theta}]_j \right\} \big](z_i)
\end{align*}
Then, we can use Lemma~\ref{lem:mp add and hom} to obtain
\begin{align*}
   \sdp Q_{\bar{\theta}}(z_i) 
   &=  \max_{j\in [1:p]} \left\{ \sdp \big[ [f]_j + [\bar{\theta}]_j \big](z_i) \right\} 
   =  \max_{j\in [1:p]} \left\{ \sdp [f]_j(z_i) + \gamma [ \bar{\theta}]_j \right\} 
   =  \max_{j\in [1:p]} \left\{ [G\sam]_{i,j} + \gamma  [ \bar{\theta}]_j \right\}.
\end{align*}
The equality $F\sam \boxtimes \bar{C}\sam = G\sam$ then implies that 
\begin{align*}
   \sdp Q_{\bar{\theta}}(z_i) 
   &=  \max_{j\in [1:p]} \left\{ [F\sam \boxtimes \bar{C}\sam]_{i,j} + \gamma  [ \bar{\theta}]_j \right\} 
   =  \max_{j\in [1:p]} \left\{ \max_{k \in [1:p]} \big\{ [F\sam]_{i,k} + [\bar{C}\sam]_{k,j} \big\} + \gamma  [ \bar{\theta}]_j \right\} \\
   & =  \max_{k\in [1:p]} \left\{ [F\sam]_{i,k} +  \max_{j \in [1:p]} \big\{  [\bar{C}\sam]_{k,j}  + \gamma  [ \bar{\theta}]_j \big\} \right\} 
   = \max_{k\in [1:p]} \left\{ [F\sam]_{i,k} +  [\bar{C}\sam \boxtimes (\gamma \bar{\theta})]_k \right\}.
\end{align*}
Finally, using the equality $\bar{\theta} = \bar{C}\sam \boxtimes (\gamma \bar{\theta})$, we have
\begin{align*}
   \sdp Q_{\bar{\theta}}(z_i) 
   &= \max_{k\in [1:p]} \left\{ [F\sam]_{i,k} +  [\bar{\theta}]_k \right\} = \max_{k\in [1:p]} \left\{ [f]_k(z_i) +  [\bar{\theta}]_k \right\} 
   = f(z_i)\tr \boxtimes \bar{\theta} 
   = Q_{\bar{\theta}}(z_i).
\end{align*}
This completes the proof.
\subsection{Proof of Proposition~\ref{prop:MP_Bellman_var}}

For each $k \in [1:q]$, we have
\begin{align*}
    \max_{i \in [1:n]} \{ [h]_k(z_i) + Q_{ \bar{\theta}}(z_i) \}  &= \max_{i \in [1:n]} \{ [h]_k(z_i) + f(z_i)\tr \boxtimes \bar{\theta}  \} \\ 
    & = \max_{i \in [1:n]} \left\{ [h]_k(z_i) + \max_{j\in [1:p]} \left\{ [f]_j(z_i) + [\bar{\theta}]_j \right\}  \right\} \\
    & = \max_{j\in [1:p]}  \left\{ \max_{i \in [1:n]} \left\{ [h]_k(z_i) +   [f]_j(z_i)\right\} + [\bar{\theta}]_j   \right\} \\
    & = \max_{j\in [1:p]}  \left\{ \max_{i \in [1:n]} \left\{ [H\sam]_{i,k} +   [F\sam]_{i,j}\right\} + [\bar{\theta}]_j   \right\} \\
    & = \max_{j\in [1:p]}  \left\{ [F\sam^{H}]_{k,j} + [\bar{\theta}]_j   \right\} = [ F\sam^{H} \boxtimes \bar{\theta}]_k.
\end{align*}
Hence, using \eqref{eq:DP_empirical_mp_var}, we can write
\begin{align*}
    \max_{i \in [1:n]} \{ [h]_k(z_i) + Q_{\bar{\theta}}(z_i) \} &= [G\sam^{H} \boxtimes (\gamma\bar{\theta})]_k
    = \max_{j\in [1:p]}  \left\{ [G\sam^{H}]_{k,j} + \gamma [\bar{\theta}]_j   \right\} \\
    & = \max_{j\in [1:p]}  \left\{ \max_{i \in [1:n]} \left\{ [H\sam]_{i,k} +   [G\sam]_{i,j}\right\} + \gamma[\bar{\theta}]_j   \right\} \\
   & = \max_{j\in [1:p]}  \left\{ \max_{i \in [1:n]} \left\{ [h]_k(z_i) +  \sdp [f]_j(z_i)\right\} + \gamma[\bar{\theta}]_j   \right\} \\
   & = \max_{i \in [1:n]}  \left\{  [h]_k(z_i) + \max_{j\in [1:p]}\left\{   \sdp [f]_j(z_i) + \gamma[\bar{\theta}]_j \right\} \right\}.
\end{align*}
Finally, we can use Lemma~\ref{lem:mp add and hom} to obtain
\begin{align*}
    \max_{i \in [1:n]} \{ [h]_k(z_i) + Q_{\bar{\theta}}(z_i) \} 
    & = \max_{i \in [1:n]}  \left\{  [h]_k(z_i) + \sdp \big[ \max_{j\in [1:p]}\left\{[f]_j + [\bar{\theta}]_j \right\}\big](z_i) \right\} \\
    & = \max_{i \in [1:n]}  \left\{  [h]_k(z_i) + \sdp Q_{\bar{\theta}}(z_i) \right\}.
\end{align*}
This completes the proof.
\subsection{Proof of Lemma~\ref{lem:mp-vqi conv}}

We first show that $F\sam^H = H\sam\tr  \boxtimes F\sam  \in \Rl^{q\times p} $ is column $\R$-astic. 
Fix $j \in [1:p]$. 
Condition~\eqref{eq:F col realistic} implies that there exists $i^* \in [1:n]$ such that $[F\sam]_{i^*,j} = [f]_j(z_{i^*}) \neq -\infty$. 
Since $H\sam$ is row $\R$-astic (by Assumption~\ref{As:compatible test}) implies that there exists $k^* \in [1:q]$ such that $[H\sam]_{i^*,k^*} = [h]_{k^*}(z_{i^*}) \neq -\infty$. 
Hence,
\begin{align*}
    [F\sam^H]_{k^*,j} &= [H\sam\tr \boxtimes F\sam ]_{k^*,j} = \max_{i \in [1:n]} \{ [H\sam]_{i,k^*} + [F\sam]_{i,j}  \} 
    \geq [F\sam]_{i^*,j} + [H\sam]_{i^*,k^*} \neq -\infty. 
\end{align*}

We next show that $F\sam^H$ is row $\R$-astic. Fix $k \in [1:q]$. 
Condition~\eqref{eq:H col realistic} implies that there exists $i^* \in [1:n]$ such that $[H\sam]_{i^*,k} = [h]_{k}(z_{i^*}) \neq -\infty$.
Also, since $F\sam$ is row $\R$-astic (by Assumption~\ref{As:features}) implies that there exists $j^* \in [1:p]$ such that $[F\sam]_{i^*, j^*} =[f]_{j^*}(z_{i^*}) > -\infty$. 
Then,
\begin{align*}
   [F\sam^H]_{k ,j^*} &=  [H\sam\tr \boxtimes F\sam]_{k, j^*} = \max_{i \in [1:n]} \{ [H\sam]_{i,k} + [F\sam]_{i, j^*} \} 
   \geq [H\sam]_{i^*,k} + [F\sam]_{i^*, j^*} \neq  -\infty.  
\end{align*}
Finally, we show that $G\sam^H = H\sam\tr \boxtimes G\sam  \in \Rl^{q\times p}  $ is row $\R$-astic. 
Fix $k \in [1:q]$. 
Condition~\eqref{eq:H col realistic} implies that there exists $i^* \in [1:n]$ such that $[H\sam]_{i^*,k} = [h]_{k}(z_{i^*}) \neq -\infty$.
Moreover, since $G\sam$ is row $\R$-astic (by Assumption~\ref{As:features}) implies that there exists $j^* \in [1:p]$ such that $[G\sam]_{i^*, j^*} = \sdp [f]_{j^*}(z_{i^*}) > -\infty$. 
Then, 
\begin{align*}
   [G\sam^H]_{k ,j^*} &=  [H\sam\tr \boxtimes G\sam]_{k, j^*} = \max_{i \in [1:n]} \{ [H\sam]_{i,k} + [G\sam]_{i, j^*} \} 
   \geq [H\sam]_{i^*,k} + [G\sam]_{i^*, j^*} \neq  -\infty. 
\end{align*}
This completes the proof.
\subsection{Proof of Theorem~\ref{thm:mp-vqi conv}} 

The proof is similar to the one provided for Theorem~\ref{thm:mp-fqi conv}.
\subsection{Proof of Theorem~\ref{thm:mp-vqi comp}} 

The compilation time includes the computation of the matrices $H\sam \in \Rl^{n\times q}$ and $F\sam,\ G\sam \in \Rl^{n\times p}$ which requires $\ord(n(p+q))$ operations. 
It also includes the computation of $F\sam^H$ and $G\sam^H$ involving MP matrix multiplications of size $q \times n \times p$ which requires $\ord(npq)$ operations. 
The compilation time complexity is then of $\ord(npq)$.
Within each iteration, the update of parameters involves a scalar multiplication of size $p$ and two MP matrix-vector multiplication of size $q \times p$ and $p \times q$.  
These computations require $\ord(pq)$ operations. 

\section{Numerical simulations}\label{app:num} 

\subsection{MDP and samples} 
The problem of interest is the stabilization of a DC motor around the origin of the state space ($x = 0$).  
The model and reward function are borrowed from~\cite[Sec.~3.4.5]{Busoniu10}. 
The system has two state variables ($d_\mathrm{x} = 2$) and one control action ($d_\mathrm{u} = 1$).
The dynamics is given by
\begin{equation}\label{eq:DC model}
x^+ = \left[ \begin{array}{cc} 1 & 0.0049 \\ 0 & 0.9540 \end{array} \right] x + \left[ \begin{array}{c} 0.0021 \\ 0.8505 \end{array} \right] u,
\end{equation}
where the state $x$ is restricted to the domain $\X = [-\pi , \pi] \times [-16\pi ,16\pi ]$` using saturation, and $u \in \U= \{ v_k = -10+5(k-1)\}_{k \in [1:5]}$ is the action space (which is, the discretization of the continuous action space~$[-10,10]$) 
The reward function is given by 
\begin{equation}\label{eq:DC reward}
r(x,u) = -x\tr \left[ \begin{array}{cc} 5 & 0 \\ 0 & 0.01 \end{array} \right] x - 0.01u^2,
\end{equation}
and the discount factor is set to $\gamma = 0.95$. 
The sample state-action pairs $\{z_i = (x_i,u_i)\}_{i \in [1:n]}$ are generated randomly from $\X\times \U$ with a uniform distribution, 
and then used to generate $x_i^+$ and $r_i$ according to dynamics~\eqref{eq:DC model} and reward function~\eqref{eq:DC reward}, respectively. 

\subsection{Features and test functions} 
For MP-FQI Algorithm~\ref{alg:mp-fqi}, the features are of the form 
$$
[f]_{jk}(x,u) = [f_{\mathrm{x}}]_j(x) +\delta_{v_k} (u), \quad j\in [1:p_{\mathrm{x}}],\ k\in [1:5].
$$ 
For the state features, we used (i)~quadratic features $[f_{\mathrm{x}}]_j(x) = -c \norm{x - y_j}_2^2$ where $ \{y_j \}_{j\in [1:p_{\mathrm{x}}]} \subset\X$ is a uniform grid with the same number of points in each state dimension, and (ii)~distance features $[f_{\mathrm{x}}]_j(x) = -c \min_{y\in \Y_j} \norm{x - y}_2^2$ where $ \{\Y_j \}_{j\in [1:p_{\mathrm{x}}]}$ is a uniform grid partitioning of~$\X$ with the same number of bins in each state dimension. 
Note that the total number of features is $p = 5 p_{\mathrm{x}}$. 
In our simulations, we set $c = (p_{\mathrm{x}})^{\frac{1}{d_{\mathrm{x}}}}$ so that the curvature of the basis functions is proportional to the number of basis functions in each dimension. 
This is a standard method in function approximation \cite{duvenaud2014kernel}. 
In general, one can also introduce a scaling factor $\alpha$ as in $c = \alpha (p_{\mathrm{x}})^{\frac{1}{d_{\mathrm{x}}}}$, where $\alpha$ is treated as a hyperparameter. 
We also note that the constant $c$ should be chosen such that it dominates the semi-convexity constant of the Q-function to ensure consistency as the number of basis functions tends to infinity~\cite{McEneaney06}. 

For v-MP-FQI Algorithm~\ref{alg:mp-fqi-var}, the features are the same as above. 
The test functions are also chosen to be the same as the features so that $q = p$.
  
For the FQI algorithm~\eqref{eq:fqi_stand_1}, the features are of the form
$$
[\Tilde{f}(x,u)]_{jk} =[\Tilde{f}_{\mathrm{x}}]_j(x) \cdot \ind{v_k} (u), \quad j\in [1:p_{\mathrm{x}}],\ k\in [1:5]. 
$$
For the state features, we used (i)~RBFs~$[\Tilde{f}_{\mathrm{x}}]_j(x) = \exp \big(-\frac{1}{c} \norm{x - y_j}_2^2 \big)$ and (ii)~indicator features~$[\Tilde{f}_{\mathrm{x}}]_j(x) =\ind{\Y_j} (x)$, where the points~$ \{y_j \}_{j\in [1:p_{\mathrm{x}}]}$, the sets~$ \{\Y_j \}_{j\in [1:p_{\mathrm{x}}]}$, and the constant~$c$ are the same as above for MP state features.

\bibliographystyle{apalike}
\bibliography{ref}

\end{document}